\documentclass[
oneside]{amsart}

\usepackage{inputenc}
\usepackage{cite}

\usepackage{amsthm}

\usepackage{stackengine}[2013-09-11]

\usepackage{amssymb}
\usepackage{amsmath}
\usepackage{amscd}
\usepackage[dvipdfmx]{graphicx}
\usepackage[dvipdfm, usenames]{color}
\usepackage[all]{xy}
\usepackage{comment}
\usepackage{setspace}
\usepackage{latexsym}
\usepackage[normalem]{ulem}
\usepackage{url}

\usepackage{xcolor}
\DeclareRobustCommand{\erase}{\bgroup\markoverwith{\textcolor{red}{\rule[.5ex]{2pt}{0.4pt}}}\ULon}

\newcommand{\mathsym}[1]{{}}
\newcommand{\unicode}[1]{{}}

\usepackage[usenames]{xcolor}

\usepackage[hyperfootnotes=false]{hyperref}

\title[]{On neighborhoods of projective space bundles over elliptic curves}
\author[]{Takayuki Koike$^{\dag}$}

\address{Department of Mathematics, Graduate School of Science, Osaka Metropolitan University, 3-3-138, Sugimoto, Sumiyoshi-ku Osaka, 558-8585 Japan}
\email{tkoike@omu.ac.jp}

\author{Laurent Stolovitch$^{\dag\dag}$}
\address{CNRS and Laboratoire J.-A. Dieudonn\'e
	U.M.R. 7351, Universit\'e C\^ote d'Azur, Parc Valrose
	06108 Nice Cedex 02, France}
\email{stolo@unice.fr}
\thanks{ This work was supported by Cooperation grant Sakura PHC 51183QC (PI: Laurent Stolovitch), JSPS Bilateral Program Number JPJSBP120243210 (PI: Takayuki Koike), by the Research Institute for Mathematical Sciences,
an International Joint Usage/Research Center located in Kyoto University, and MEXT Promotion of Distinctive Joint Research Center Program JPMXP0723833165 (Osaka Central Advanced Mathematical Institute, Osaka Metropolitan University).
The first author is supported by the Grant-in-Aid for Scientific Research C (23K03119) from JSPS. The authors thank Dr. Xiaojun Wu for helpful and insightful discussions.}

\newcommand{\diag}{\operatorname{diag}}

\theoremstyle{plain} 
\newtheorem{theorem}{\noindent\bf Theorem}[section]
\newtheorem{lemma}[theorem]{\noindent\bf Lemma}

\newtheorem{proposition}[theorem]{\noindent\bf Proposition}

\theoremstyle{definition} 
\newtheorem{definition}[theorem]{\noindent\bf Definition}
\newtheorem{remark}[theorem]{\noindent\bf Remark}

\newcommand{\ga}{\begin{gather}}
	\newcommand{\ega}{\end{gather}}
\newcommand{\gan}{\begin{gather*}}
	\newcommand{\egan}{\end{gather*}}
\newcommand{\al}{\begin{align}}
	\newcommand{\eal}{\end{align}}
\newcommand{\aln}{\begin{align*}}
	\newcommand{\ealn}{\end{align*}}
\newcommand{\eq}[1]{\begin{equation}\label{#1}}
	\newcommand{\eeq}{\end{equation}}

	\def\beq{\begin{equation}}
		\def\eeq{\end{equation}}

\newcommand{\N}{\mathbb{N}}
\newcommand{\C}{{\mathbb C}}

\newcommand{\Z}{{\mathbb Z}}

\newcommand{\cyl}{{\widetilde C}}

\newcommand{\IM}{\operatorname{Im}}

\newcommand{\cL}{\mathcal}

\newcommand{\del}{\delta}
\newcommand{\Del}{\Delta}

\newcommand{\e}{\epsilon}
\newcommand{\om}{\omega}
\newcommand{\Om}{\Omega}

\newcommand{\la}{\lambda}

\newcommand{\re}[1]{(\ref{#1})}

\newcommand{\rp}[1]{Proposition~\ref{#1}}

\begin{document}

\begin{abstract}
We give conditions ensuring that a neighborhood of an embedded projective space bundle over an elliptic curve is holomorphically equivalent to a neighborhood of the zero section of its normal bundle. 
\end{abstract}

\date{\today}
 \maketitle



	\section{Introduction}
Given a compact complex manifold $C$ embedded into a complex manifold $M$, one obtains another embedding $C$ as the zero section of the (holomorphic) normal bundle $N_{C/M}$ of $C$ into $M$. A natural problem is to compare the neighborhoods of $C$ in both manifolds. Following somehow Grauert's Formale Prinzip, we are interested to understand under which conditions such two neighborhoods are first of all, {\it formally equivalent} and if so, {\it holomorphically equivalent}. This problem has definitive answers whenever the normal bundle $N_{C/M}$ admits a Hermitian fiber metric whose curvature is sufficiently positive or negative in a suitable sense (e.g \cite{grauert-neg,hironaka-rossi,G}). In this article we shall focus on the case where the normal bundle is {\it flat}. Since the pioneering work of Arnol'd \cite{arnold-elliptic, ilyashenko-tori},  there has been very few other results related to this situation (e.g \cite{GS21,stolo-gong-tori}). Another related problem is the so-called ``vertical linearization" which ensure the existence of a holomorphic foliation in a neighborhood of $C$, having $C$ as a leaf. Following the pioneering work of Ueda \cite{U}, considering curve embedded into a surface, there has been some works in order to extend the result to either higher dimensional or co-dimensional $C$ (e.g. \cite{Koi,GS21,stolo-wu}).
	
	 The aim of this article is built another completely new example of such embeddings using previously known results of the kind. Let us give our framework~:

\begin{itemize}
\item $C$ is a compact complex torus of dimension $n$ (we consider mainly in the case where $n=1$). 
\item $Y$ is the total space of a $\mathbb{P}^r$-bundle $\pi\colon Y\to C$ over $C$. 
\item $Y$ is holomorphically embedded into a complex manifold $X$ . 
\end{itemize}
Note that, if $n=1$, there exists a vector bundle $E\to C$ of rank $r+1$ such that $Y$ coincides with the fiberwise projectivization $\mathbb{P}(E)$ of $E$ (i.e. the quotient of the complement of the zero section in the total space of $E$ by the fiberwise Hopf action~: {$x \in C$, ${\mathbf e}_x\in E_x$, $\lambda\in \C^*\mapsto \lambda{\mathbf e}_x$). 
The existence of such $E$ can be shown by the same argument in the proof of Proposition \ref{prop:star_Pr_bdl_over_ellipt_curve} below. 
\begin{center}
\includegraphics[width=6cm]{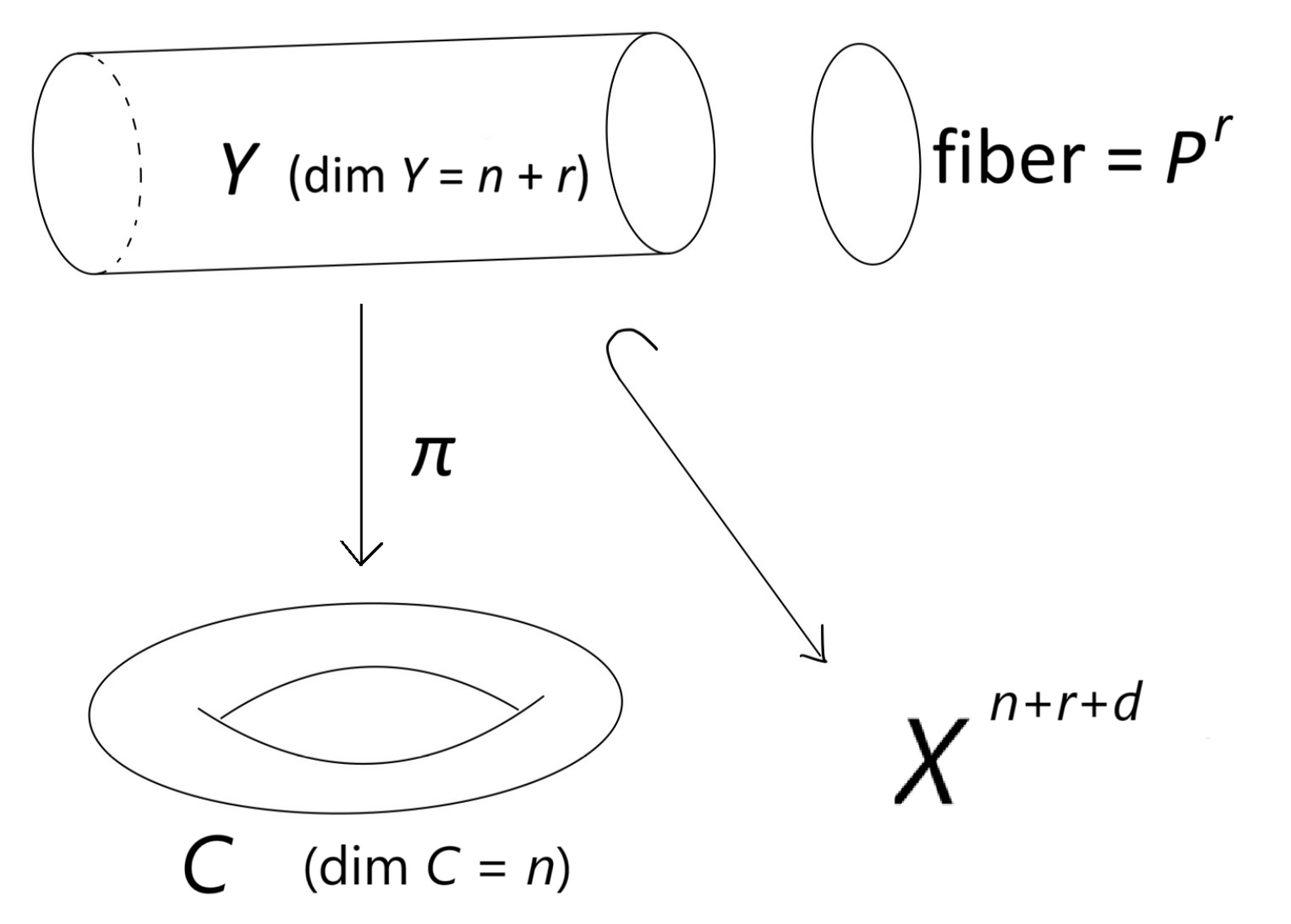}
\end{center}

Note that $\pi_1(Y, *)$ is isomorphic to $\pi_1(C, *)\cong \Z^{2n}$, 
from which it follows that any unitary flat vector bundle $E$ of rank $R=:r+1$ over $Y$ is isomorphic to the pull-back $\pi^*E_C$ for some unitary flat vector bundle $E_C$ over $C$ of the same rank (See Lemma \ref{lem:simplyconnfibration_fund_group}). Therefore, such a vector bundle $E_C$ is completely determined by $n$ diagonal unitary $R\times R$-matrices $M_i=\diag(\mu_{i,1},\ldots,\mu_{i,R} )$, $i=1\ldots, n$ (see e.g. \cite{iena-05}).
In what follows, we shall mainly focus on the case of an elliptic curve $C\simeq  \C/\langle 1,\tau\rangle$ where $\IM\tau\neq 0$. Let us set $\lambda:=e^{2i\pi\tau}$. We will say that a unitary flat vector bundle $E$ over $Y$ is {\it Diophantine} if the corresponding vector bundle $E_C$ satisfies the (non-resonant) Diophantine condition ~: 
\begin{definition}[\cite{stolo-gong-tori}]
A flat unitary vector bundle $E_C$ over $C$ is said to be (resp. strongly) vertically Diophantine if, for all $p\in\Z$, $Q\in \N^d$,  $\ell=1,\ldots, d$
\beq\label{vd}
\max_i|\lambda^p\mu_{i}^Q-\mu_{i,\ell}|>\frac{C}{|p|+|Q|}.
\eeq
(resp.
\beq\label{svd}
\forall i=1,\ldots, n,\quad |\lambda^p\mu_{i}^Q-\mu_{i,\ell}|>\frac{C}{|p|+|Q|}.
\eeq)\\
The vector bundle $E_C$ over $C$ is said to be (non-resonant) Diophantine if in addition to \re{vd} it also satisfies to 
\begin{align}\label{hd}
\max_i|\lambda^p\mu_{i}^Q-\lambda|&>\frac{C}{|p|+|Q|}.
\end{align}
\end{definition}
This definition has been devised in the case of complex torus in any dimension in \cite{stolo-gong-tori}.
\begin{remark}
	As $\IM\tau\neq 0$ and $|\mu_{i,\ell}|=1$ for all $i=1,\ldots, n$, $\ell=1,\ldots, R$, conditions \re{vd} and \re{hd}, barely boil down to the conditions with $p=0$. Indeed, otherwise $\lambda^p\mu_{i}^Q$ are never close neigher to $\mu_{i,\ell}$ nor to $\lambda$.
\end{remark}
\begin{remark}
	In the case of an elliptic curve (i.e. $n=1$), the strongly vertically Diophantine condition is identical to vertically Diophantine condition.
\end{remark}
 Also, we say that a unitary flat vector bundle $E$ over $Y$ is {\it vertically Diophantine} if the corresponding vector bundle $E_C$ satisfies vertically Diophantine.
 \begin{definition}
 	A neighborhood of $C$ in $M$ is said to be  holomorphically full-linearizable if it is holomorphically equivalent to a neighborhood of the zero section in the normal bundle $N_{C/M}$. 
 	\end{definition}


Our main result in the present paper is the following: 
\begin{theorem}\label{cor_new}
Let $Y$ be the direct product of an elliptic curve and the projective space. 
Assume that $Y$ is embedded in a complex manifold $X$. 
Assume also that the (holomorphic) normal bundle $N_{Y/X}$ is a Diophantine unitary flat vector bundle, and that the restriction $T_X|_Y$ of the holomorphic tangent bundle $T_X$ splits as $T_Y\oplus N_{Y/X}$. 
Then the neighborhood of $Y$ in $X$ can be full-linearized. 
\end{theorem}

Theorem \ref{cor_new} follows from the following more general theorem, 
which is a generalization of the linearizability theorems for a neighborhood of tori \cite{A} \cite{stolo-gong-tori}. 
\begin{theorem}\label{cor:main}
Assume that $n=1$ (i.e. $C$ is an elliptic curve). Let $Y$ be a projective bundle over $C$ that coincides with the fiberwise projectivization $\mathbb{P}(E)$ of a vector bundle $E\to C$ of rank $r+1$. Assume that $Y$ is embedded in a complex manifold $X$. 
Assume also that the (holomorphic) normal bundle $N_{Y/X}$ is a Diophantine unitary flat vector bundle, the restriction $T_X|_Y$ of the holomorphic tangent bundle $T_X$ splits as $T_Y\oplus N_{Y/X}$, and that the triple $(C, E, N)$ satisfies the following condition, where $N\to C$ is the unitary flat vector bundle whose pull-back is isomorphic to $N_{Y/X}\to Y$.  
\begin{description}
\item[(Condition $\ast$)] For any open neighborhood $M$ of the zero section in the total space of $N$ and any holomorphic vector bundle $F\to M$ of rank $r+1$ whose restriction to the zero section is holomorphically isomorphic to $E$ (via the natural identification of $C$ with the zero section), there exists a smaller neighborhood $M'\subset N$ such that the restriction $F|_{M'}$ is holomorphically isomorphic to the pull-back $(P|_{M'})^*E$, where $P\colon N\to C$ is the projection. 
\end{description}
\noindent Then the neighborhood of $Y$ in $X$ can be full-linearized. 
\end{theorem}

The strategy of the proof of Theorem \ref{cor:main} is as follows. 
First we observe a neighborhood of each fibers of $\pi$ in $X$. 
By applying a higher codimensional analogue \cite[Theorem 1.1]{Koi} of Ueda's theorem \cite{U} and the well-known deformation rigidity property of projective spaces (see Lemma \ref{lem:P^r_satisfies_assumption123} below), we construct nice atlas on a neighborhood of $Y$ such that the transitions never depend on the coordinates in the $\pi$-fiber direction (\S \ref{subsection:constr_of_nice_atlas_on_a_nbhd_Y}). 
By using the transitions of this atlas, we construct a manifold $M$ of dimension $n+d$ which includes $C$ as a submanifold, together with a $\mathbb{P}^r$-fibration $V\to M$ from a neighborhood $V$ of $Y$ in $X$ whose restriction to $Y$ coincides with $\pi$ (\S \ref{subsection:constr_of_CM}, ``smashing'' of the $\mathbb{P}^r$'s in $V$). 
The condition $n=1$ is needed only for showing that this $\mathbb{P}^r$-fibration is a fiberwise projectivization of a vector bundle on $M$. 
By applying \cite[Theorem 1.1]{stolo-gong-tori}, there is a neighborhood of $C$ in $M$ holomorphically equivalent to a neighborhood of $C$ in its normal bundle. 
We shall apply {\bf (Condition $\ast$)} in \rp{prop:star_Pr_bdl_over_ellipt_curve}, to finally linearize the neighborhood of $Y$ in $X$. 

To show Theorem \ref{cor_new}, we shall prove in section \ref{section:prf_of_thm_extension-trivial-tori}, the following theorem~:

\begin{theorem}\label{thm:extension-trivial-tori}
Let $N_C$ be a unitary flat vector bundle over a compact complex torus $C$, 
and $M$ be a neighborhood of the zero section in the total space of $N_C$. 
Let ${\mathcal E}\rightarrow M$ be a holomorphic vector bundle of rank $\ell$ such that ${\mathcal E}|_{C}$ is holomorphically trivial, where the zero section is identified with $C$. 
Assume that $N_C$ is vertically Diophantine. 
Then there exists a neighborhood $M'$ of $C$ in $M$ such that ${\mathcal E}|_{M'}\rightarrow M'$ is holomorphically trivial.
\end{theorem}
We shall apply it to show that 
{\bf (Condition $\ast$)} is satisfied when $Y$ is the direct product of an elliptic curve and the projective space. 

Note that Theorem \ref{thm:extension-trivial-tori} is a generalization of \cite[Proposition 3.3]{KU} and \cite[Lemma 4.6]{Koi2}. 
Theorem \ref{thm:extension-trivial-tori} is shown by applying Newton scheme to construct a nice local trivializations of a given vector bundle. 

%
\section{Some fundamentals}

\subsection{On the fundamental group of a fiber bundle with simply-connected fibers}

\begin{lemma}\label{lem:simplyconnfibration_fund_group}
Let $Z$ be a connected and simply-connected manifold, and $Y$ be the total space of a $Z$-fiber bundle over a connected manifold $C$. 
Then it holds that $\pi_1(Y, *)\cong \pi_1(C, *)$. 
Moreover, any unitary flat vector bundle on $Y$ is isomorphic to the pull-back of a unitary flat vector bundle on $C$.  
\end{lemma}

\begin{proof}
It can be shown by using the exactness of the sequence 
\[
\cdots \to \pi_2(C, *)\to \pi_1(Z, *)\to\pi_1(Y, *) \to \pi_1(C, *)\to \pi_0(Z, *) \to \cdots.
\]
The latter half of the assertion follows from the correspondence between the unitary flat vector bundles and the unitary representations of the fundamental groups (See \cite[Proposition 2.2]{Koi} for example). 
\end{proof}

\subsection{Properness on a neighborhood of certain submanifold}
Here we show the following: 
\begin{lemma}\label{lem:properness_fund}
Let $A$ and $B$ be differentiable manifolds, 
$\Xi\subset A$ a compact differentiable submanifold, and $f\colon A\to B$ be a continuous map such that $\Xi=f^{-1}(p)$ holds for a point $p\in B$. 
Then there exists a coordinate open ball of $B$ centered at $p$ such that, for the connected component $V$ of the preimage $f^{-1}(B)$ which includes $\Xi$, the map $f|_V\colon V \to B$ is proper. 
\end{lemma}

\begin{proof}
Take a ($C^\infty$-) tubular neighborhood $T$ of $\Xi$ in $A$. 
As the boundary $\partial T$ is homeomorphic to a sphere bundle over $\Xi$, $\partial T$ is compact. Thus the image $K:=f(\partial T)$ is compact. 
We have $p\not\in K$, since $f^{-1}(p)=\Xi$ does not intersect $\partial T$. 
As $B$ is metrizable, $B$ is a regular Hausdorff space (i.e. $B$ is a $T_3$-space). Therefore there exists a coordinate open ball $B'$ of $B$ centered at $p$ such that $B'\cap K=\emptyset$. As $f$ is continuous and $f^{-1}(K')$ is closed for any compact subset $K'\subset B'$, we have that $f^{-1}(K')\cap\overline{T}$ is a closed subset of $\overline{T}$ which does not intersects $\partial T$. 
Thus $f|_V\colon V\to B'$ is proper, where $V$ is the connected component of $f^{-1}(B')$ which includes $\Xi$. 
\end{proof}


\subsection{Some properties of projective space bundles}\label{subsection:someFanoBdls}
In this subsection, we let $C$ be a complex manifold of dimension $n$ and $\pi\colon Y\to C$ be a holomorphic fiber bundle over $C$ whose fibers are biholomorphic to a compact complex manifold $Z$ of dimension $r$ which satisfies the following three assumptions: 
\begin{description}
  \setlength{\parskip}{0cm} 
  \setlength{\itemsep}{0cm} 
\item[(Assumption 1)] $Z$ is simply connected. 
\item[(Assumption 2)] $Z$ is rigid in the sense of deformation: i.e. for any deformation (proper holomorphic surjection) $p\colon \mathcal{Z}\to S$ between complex manifolds and a point $0\in S$ whose fiber $Z_0:=p^{-1}(0)$ is biholomorphic to $Z$, there exists a neighborhood $S'$ of $0$ in $S$ and a biholomorphism $\mathcal{Z}|_{S'}:=p^{-1}(S')\cong Z \times S'$ which makes the diagram 
\[
\xymatrix{
\mathcal{Z}|_{S'} \ar[rd]_{p|_{\mathcal{Z}}|_{S'}} \ar[rr]^{\cong} &  & Z\times S' \ar[ld]^{{\rm Pr}_2} \\
& S' & 
}
\]
commutative, where ${\rm Pr}_2$ is the second projection. 
\item[(Assumption 3)] It holds that $H^1(Z, \mathcal{O}_Z)=0$. 
\end{description}
One can apply the facts shown in this subsection to a projective space bundle over $C$ by virtue of the following: 

\begin{lemma}\label{lem:P^r_satisfies_assumption123}
$Z=\mathbb{P}^r$ satisfies the three assumptions above. 
\end{lemma}

\begin{proof}
As $\mathbb{P}^r$ admits a cell decomposition which consists of a single cell in each even degree, it is simply connected. 
From \cite[Theorem 5.1]{KO}, it follows that $H^q(\mathbb{P}^r, T_{\mathbb{P}^r})=0$ for $q=1, 2$ and $H^1(\mathbb{P}^r, \mathcal{O}_{\mathbb{P}^r})=0$. 
By applying the vanishing of $H^q(\mathbb{P}^r, T_{\mathbb{P}^r})$ for $q=1, 2$, it follows from Kodaira--Nirenberg--Spencer theorem \cite{KNS} (see also {\cite[p. 45 Theorem 3.2]{KM}} or {\cite[Theorem 2.5, p. 109]{H}}) that there exists a versal deformation such that the base space is a (non-singular) manifold of dimension $0$. Thus the rigidity holds in the sense of deformation. 
\end{proof}

For $\pi\colon Y\to C$ as above, first we show the following:
\begin{lemma}\label{lem:unitary_flat_triv_on_Z}
Any unitary flat vector bundle on $Z$ is holomorphically trivial. 
\end{lemma}

\begin{proof}
Let $E$ be a unitary flat vector bundle on $Z$ of rank $d$. 
By the correspondence between the unitary flat vector bundles and the unitary representations of the fundamental groups (See \cite[Proposition 2.2]{Koi} for example), it follows that there exists a group representation $\sigma'\colon \pi_1(Z, *)\to {\rm U}(d)$ such that 
\[
E\cong \widetilde{Z}\times\mathbb{C}^d/\sim_{\sigma'}
\]
holds, where $\widetilde{Z}$ is the universal covering of $Z$ and ``$\sim_{\sigma'}$'' is the relation such that $(z, v)\sim_{\sigma'} (z', v')$ holds iff $z'=\gamma(z)$ and $v'=\sigma'(\gamma)\cdot v$ holds for some deck transformation $\gamma$. 
From {\bf (Assumption 1)}, it follows that $\sigma'$ is the trivial representation, from which it follows that $E\cong \mathcal{O}_Z^{\oplus d}$. 
\end{proof}

Next let us show the following:
\begin{lemma}\label{lem:normal_bdl_triv_on_Z}
Assume that $Y$ is as above, and that $Y$ is embedded in a complex manifold $X$ so that the normal bundle $N_{Y/X}$ is a unitary flat vector bundle on $Y$ of rank $d$. Then, for a point $p\in C$, the normal bundle $N_{Y_p/X}$ of the fiber $Y_p:=\pi^{-1}(p)$ in $X$ is holomorphically trivial: i.e. $N_{Y_p/X}\cong \mathcal{O}_{Y_p}^{\otimes n+d}$. 
\end{lemma}

\begin{proof}
Consider the short exact sequence 
\[
0\to N_{Y_p/Y} \to N_{Y_p/X} \to N_{Y/X}|_{Y_p}\to 0. 
\]
\begin{figure}[h]
\begin{center}
\includegraphics[width=6cm]{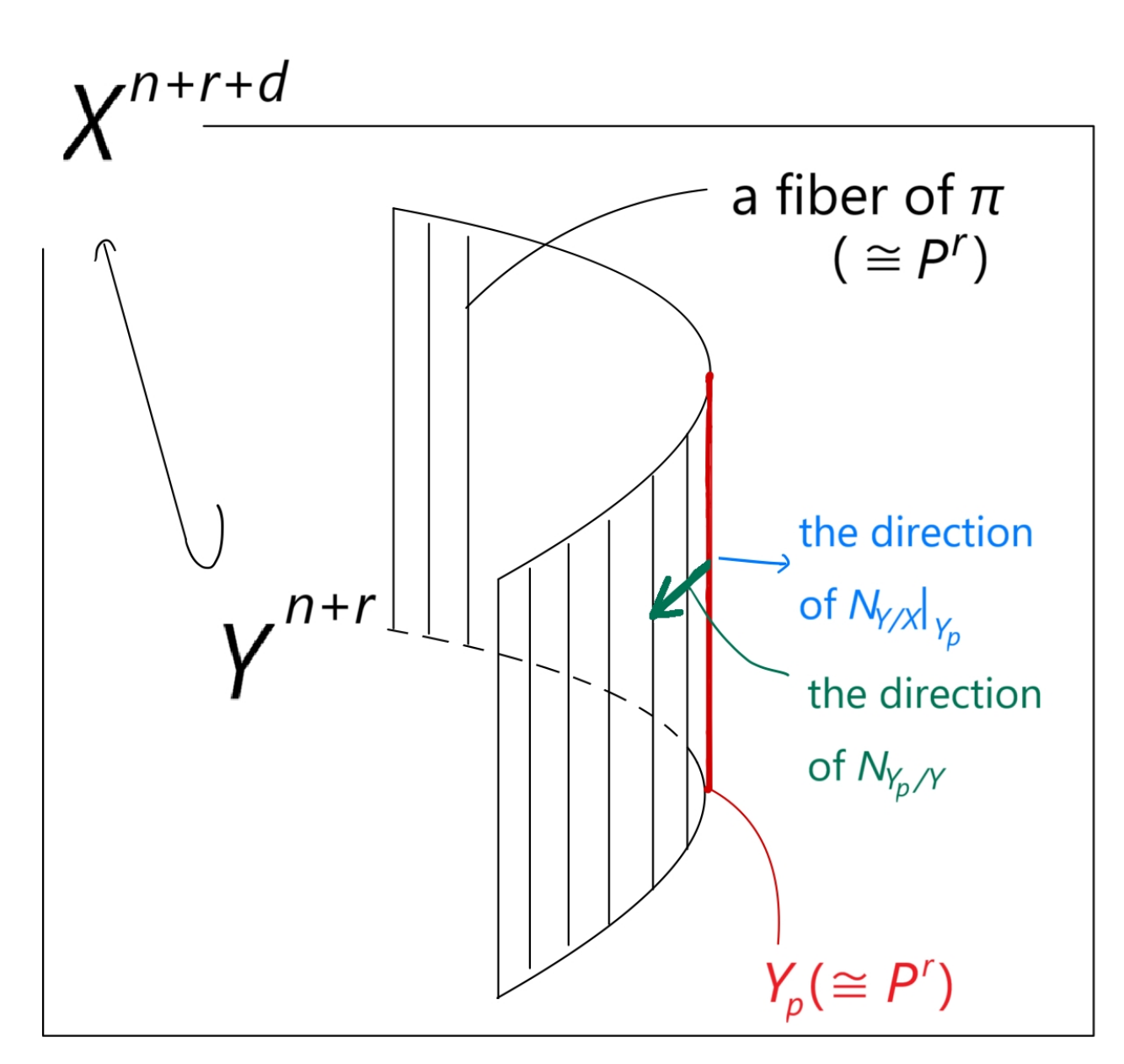}
\end{center}
\end{figure}
 As $\mathbb{P}^r$ satisfies {\bf (Assumption 2)} by Lemma \ref{lem:P^r_satisfies_assumption123}, 
there exists a neighborhood $S'$ of $p$ in $C$ and a biholomorphism $W:=\pi^{-1}(S')\cong \mathbb{P}^r\times S'$ which makes the diagram 
\[
\xymatrix{
W \ar[rd]_{\pi|_{W}} \ar[rr]^{\cong} &  & \mathbb{P}^r \times S'\ar[ld]^{{\rm Pr}_2} \\
& S' & 
}
\]
commutative, where ${\rm Pr}_2$ is the second projection. 
As $Y_p$ corresponds to the fiber $\{p\}\times \mathbb{P}^r$, we have 
$N_{Y_p/Y}\cong \mathcal{O}_{Y_p}^{\oplus n}$. 
From Lemma \ref{lem:unitary_flat_triv_on_Z}, it follows that $N_{Y/X}|_{Y_p}\cong \mathcal{O}_{Y_p}^{\oplus d}$. As the extension class of the short exact sequence above is an element of 
\[
{\rm Ext}^1(N_{Y/X}|_{Y_p}, N_{Y_p/Y})\cong 
H^1(Y_p, (N_{Y/X}|_{Y_p})^*\otimes N_{Y_p/Y})
\cong H^1(Z, \mathcal{O}_Z^{\oplus nd})=H^1(Z, \mathcal{O}_Z)^{\oplus nd}, 
\]
which is zero by {\bf (Assumption 3)}. 
The assertion follows from this. 
\end{proof}

\section{Proof of Theorem \ref{cor:main}}\label{section:prf_of_thmmain}
In this section, we show Theorem \ref{cor:main}. 

\subsection{Construction of nice atlas on a neighborhood of $Y$}\label{subsection:constr_of_nice_atlas_on_a_nbhd_Y}

First let us show the following: 
\begin{lemma}\label{lem:full-linearizability_Z}
Assume that $\pi\colon Y\to C$ is as in \S \ref{subsection:someFanoBdls} with $Z=\mathbb{P}^r$, and that $Y$ is embedded in a complex manifold $X$ so that the normal bundle $N_{Y/X}$ is a unitary flat vector bundle on $Y$ of rank $d$. Then, for a point $p\in C$ and its fiber $Y_p:=\pi^{-1}(p)$, there exits an open neighborhood $V_p$ of $Y_p$ in $X$ and a biholomorphism $F_p\colon V_p\cong Z\times \Delta^{n+d}$ such that $F_p|_{Y_p}\colon Y_p\cong Z\times \{0\}$ holds. 
\end{lemma}

\begin{center}
\includegraphics[width=6cm]{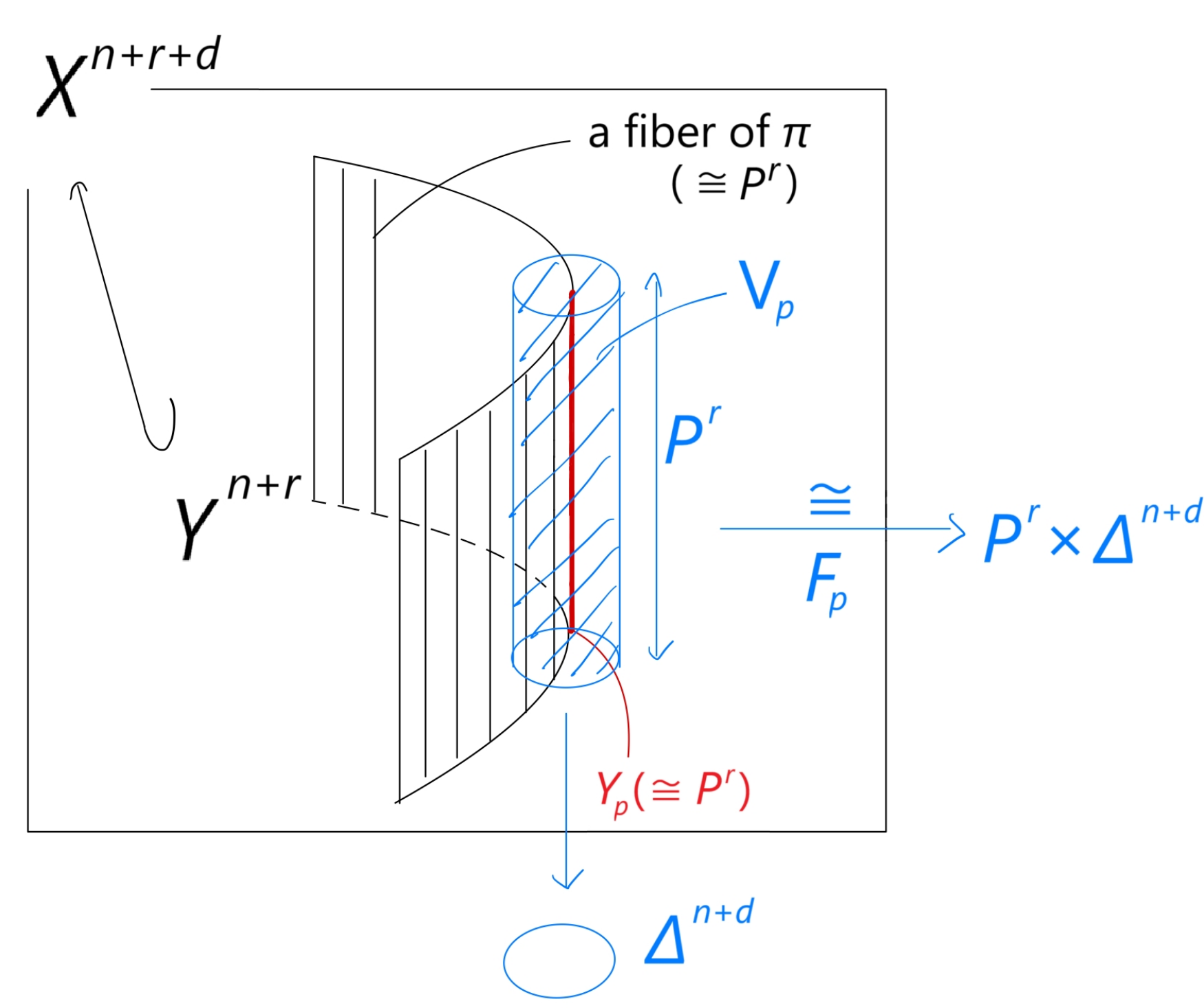}
\end{center}

\begin{proof}
From Lemma \ref{lem:normal_bdl_triv_on_Z} and {\bf (Assumption 3)}, it follows that 
\[
H^1(Y_p, N_{Y_p/X}\otimes S^\ell (N_{Y_p/X})^*)=0
\]
holds for any positive integer $\ell$. 
Therefore we have that $(Y_p, X)$ is of infinite type in the sense of \cite{Koi}. 
By a higher (co-)dimensional analogue of Ueda's theorem (see e.g. \cite{Koi}[Theorem 1.1], \cite{stolo-wu}), we have that neighborhoods of $Y_p$ in $X$ can be (not only formally, but also holomorphically) vertically linearizable, or equivalently, that there exists a neighborhood $V_p$ of $Y_p$ in $X$, a open covering $\{V_{pj}\}$ of $V_p$, and a system of defining holomorphic functions $w_{pj}=(w_{pj}^1, w_{pj}^2, \dots, w_{pj}^{n+d})$ such that $w_j=T_{jk}w_k$ holds on each $V_{pj}\cap V_{pk}$, where $T_{jk}$'s are unitary matrix of order $d$ such that $N_{Y_p/X}^*=[\{(V_{pj}\cap V_{pk}, T_{jk})\}]$. 
As $N_{Y_p/X}$ is holomorphically trivial, we may assume that $T_{jk}$ is the identity matrix. Therefore we may assume by changing the order of the suffixes that $w_{pj}^\nu=w_{pk}^\nu$ holds on each $V_{jk}$ for any $\nu=1, 2, \dots, n+d$, which implies that the map (fibration) $\pi_p\colon V\to \mathbb{C}^{n+d}$ defined by 
\[
\pi_p(x) := (w_{pj}^1(x), w_{pj}^2(x), \dots, w_{pj}^{n+d}(x))\in \mathbb{C}^{n+d}\quad (\text{if}\ x\in V_{pj})
\]
is well-defined as a holomorphic function. 
By shrinking $V_p$ to the preimage of a small neighborhood of $0$ in $\mathbb{C}^{n+d}$ and by changing the scaling, by virtue of Lemma \ref{lem:properness_fund}, we may assume that $\pi_p(V)=\Delta^{n+d}$ and that $\pi_p$ is proper and surjective, where $\Delta$ denotes the unit disc in $\mathbb{C}$. As is clear by the construction, $\pi_p$ is holomorphic submersion. 
Thus the full-linearizablity of a neighborhood of $Y$ follows from Lemma \ref{lem:normal_bdl_triv_on_Z} and {\bf (Assumption 2)}. 
\end{proof}

Next we show the following:
\begin{lemma}\label{lem:nice_atlas}
Assume that $Y$ is as above, and that $Y$ is embedded in a complex manifold $X$ so that the normal bundle $N_{Y/X}$ is a unitary flat vector bundle on $Y$ of rank $d$. Then, for a point $p\in C$ and its fiber $Y_p:=\pi^{-1}(p)$, by shrinking $V_p$, there exists a neighborhood $U_p$ of $p$ in $C$ and a decomposition $\Delta^{n+d}=\Delta^n\times \Delta^d$ such that 
\[
\xymatrix{
V_p \ar[r]^{\cong\ \ \ \ \ \ }_{F_p\ \ \ \ \ \ }& Z\times \Delta^n \times \Delta^d\\
V_p\cap Y=\pi^{-1}(U_p)\ar[r]^{\cong}_{\ \ F_p|_{V_p\cap Y}\ }\ar[u]_{\rm inclusion}\ar[d]_{\pi}& Z\times \Delta^n \times \{0\}\ar[u]_{\rm inclusion}\ar[d]_{{\rm Pr}_2}\\
U_p\ar[r]_{\cong} & \Delta^n
}
\]
is commutative, where $F_p$ is the map as in Lemma \ref{lem:full-linearizability_Z}. 
\end{lemma}

\begin{center}
\includegraphics[width=6cm]{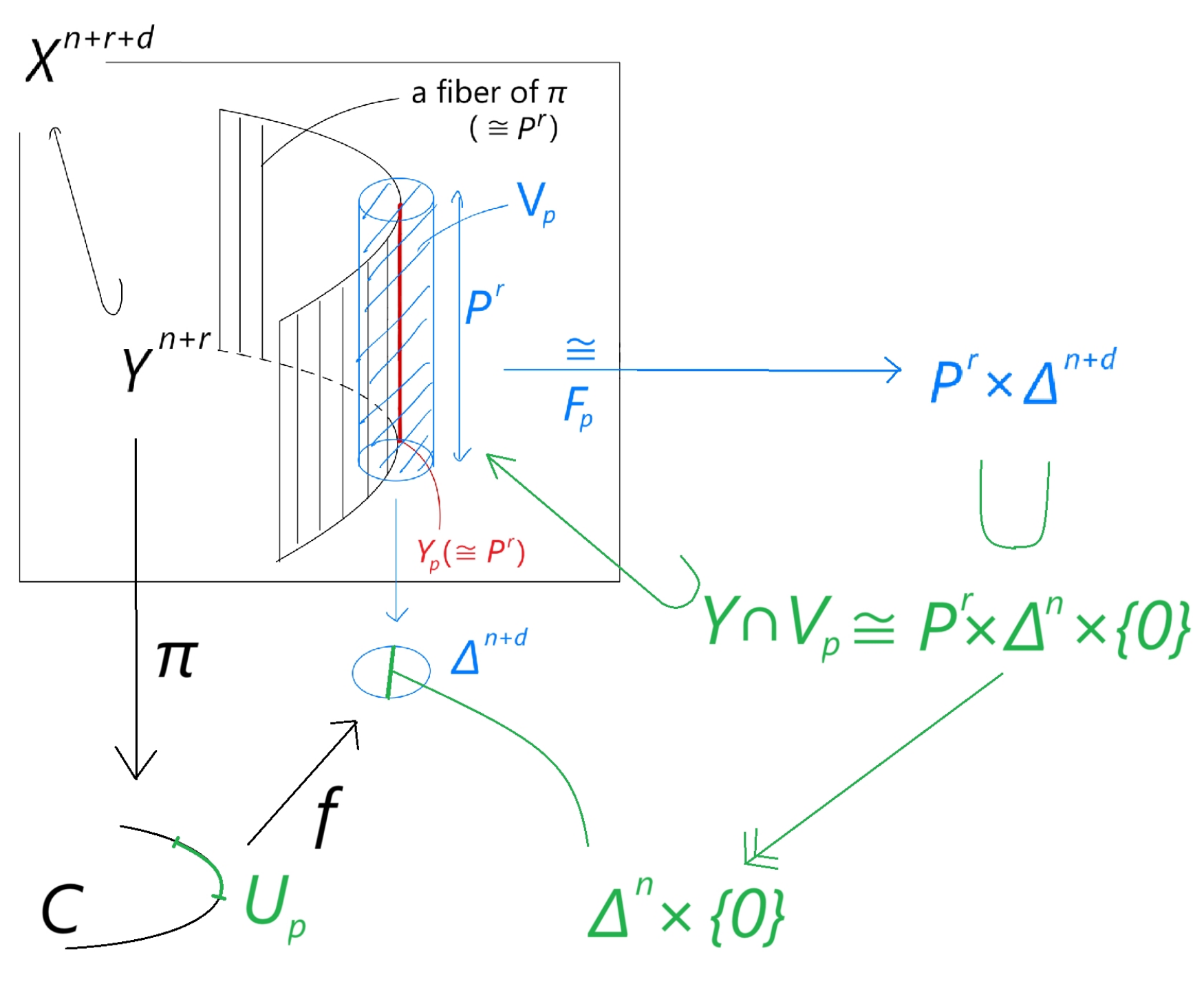}
\end{center}

\begin{proof}
From Lemma \ref{lem:properness_fund}, it follows that there exists a neighborhood $U_p$ of $p$ in $U$ such that $\pi^{-1}(U_p)\subset V_p$. Let us show that there exists a holomorhic map $f\colon U_p\to \Delta^{n+d}$ such that $f\circ \pi={\rm Pr}_2\circ F_p$ holds on $\pi^{-1}(U_p)$, where ${\rm Pr}_2\colon Z\times \Delta^{n+d}\to \Delta^{n+d}$ is the second projection: 
\[
\xymatrix{
V_p \ar[r]^{\cong\ \ \ \ \ \ }_{F_p\ \ \ \ \ \ }& Z\times \Delta^{n+d}\ar[r]^{\ \ {\rm Pr}_2} & \Delta^{n+d}\\
V_p\cap Y=\pi^{-1}(U_p)\ar[u]_{\rm inclusion}\ar[d]_{\pi}&  & \\
U_p\ar[rruu]_{\exists f} &  & 
}
\]
In order to show this, it is sufficient to show that ${\rm Pr}_2\circ F_p|_{\pi^{-1}(q)}$ is the constant function for each $q\in U_p$ (since then we can define $f$ by letting $f(q)$ be the value of ${\rm Pr}_2\circ F_p|_{\pi^{-1}(q)}$). As the constantness of the function ${\rm Pr}_2\circ F_p|_{\pi^{-1}(q)}$ is clear from the compactness of $\pi^{-1}(q)\cong Z$ and the maximal principle, the assertion holds. 

Next we will check the smoothness of the image of $f$. 
For that purpose, we will calculate the rank of the differential $(Df)_p$ of the map $f$ at the point $p$. Note that the diagram 
\[
\xymatrix{
T_{V_p, x} \ar[r]^{\cong\ \ \ \ \ \ }_{(DF_p)_x\ \ \ \ \ \ }& T_{Z\times \Delta^{n+d}, F_p(x)}\ar[r]^{\rm \ \ surj}&T_{\Delta^{n+d}, 0}\\
T_{\pi^{-1}(U_p), x}\ar[d]^{(D\pi)_x\ \ }\ar[u]_{\rm inclusion}& & \\
T_{U_p, p}\ar[rruu]_{(Df)_p} & & 
}
\]
is commutative for each point $x\in Y_p$. 
 As the map $T_{V_p, x} \to T_{\Delta^{n+d}, 0}$ is surjective, this map is of rank $n+d$. As $T_{\pi^{-1}(U_p), x}$ is a subspace of $T_{V_p, x}$ of codimension $d$, the rank of the map $T_{\pi^{-1}(U_p), x}\to T_{\Delta^{n+d}, 0}$ is larger than or equal to $(n+d)-d=n$. As ${\rm dim}\,T_{U_p, p}=n$, one has ${\rm rank}\,(Df)_p=n$. 

Therefore, by a coordinate change of $\Delta^{n+d}$, we may assume that 
\[
\text{Image}(Df)_p=\left\langle \frac{\partial}{\partial w_1},\ \frac{\partial}{\partial w_2},\ \dots,\ \frac{\partial}{\partial w_n} \right\rangle, 
\]
since ${\rm dim}\,U_p=n$. 
By virtue of Lemma \ref{lem:properness_fund}, we may assume that $f$ is proper. Thus it follows from Remmert's proper mapping theorem (see \cite[p.118 (8.8)]{D} for example) that $f(U_p)$ is an analytic space holomorphically included in $\Delta^{n+d}$. 
Moreover, as $(Df)_p$ is of full-rank, we have that $f(U_p)$ is a submanifold of dimension $n$ which is biholomorphic to $U_p$ by shrinking $U_p$ if necessary. 

Therefore, by shrinking $U_p$ so that $U_p\cong \Delta^n$ and by changing the coordinates of $\Delta^{n+d}$ suitably, we may assume that $f$ coincides with the natural inclusion $\Delta^n\to \Delta^{n+d}$ and the image of $f$ coincides with $\Delta^n\times\{(0, 0, \dots, 0)\}\subset \Delta^{n+d}$. Thus, by shrinking $V_p$ to a small neighborhood of $\pi^{-1}(U_p)$, it holds that $V_p\cap \pi^{-1}(U_p)=F_p^{-1}(Z\times \Delta^n\times\{(0, 0, \dots, 0)\})$, from which the lemma holds. 
\end{proof}

Finally we obtain the following: 
\begin{lemma}\label{lem:nice_atlas_coord_change}
Assume that $Y$ is as above, and that $Y$ is embedded in a complex manifold $X$ so that the normal bundle $N_{Y/X}$ is a unitary flat vector bundle on $Y$ of rank $d$. 
Take a point $p\in C$. Let $U_p, V_p, F_p\colon V_p\to Z\times U_p\times \Delta^{d}$ be those as in Lemma \ref{lem:nice_atlas}. We will denote by $h_p=(h_p^1, h_p^2, \dots, h_p^n)$ the coordinates of $U_p\cong \Delta^n$, and by $v_p=(v_p^1, v_p^2, \dots, v_p^d)$ the coordinates of $\Delta^d$. 
Take another point $q\in C$, and take $U_q, V_q, F_q\colon V_q\to Z\times U_q\times \Delta^{d}, h_q$, and $v_q$ in the same manner. Assume that $U_p\cap U_q\not=\emptyset$ and $V_p\cap V_q$ is a connected neighborhood of $\pi^{-1}(U_p\cap U_q)$ such that $V_p\cap V_q\cap Y=\pi^{-1}(U_p\cap U_q)$. 
Then the transition functions $\Phi^h_{qp}$ and $\Phi^v_{qp}$ from $V_p$ to $V_q$ defined by 
\[
h_q=\Phi^h_{qp}(h_p, v_p)\quad\text{and}\quad v_q=\Phi^v_{qp}(h_p, v_p)
\]
do not depend on the coordinate of $Z$. 
\end{lemma}

\begin{proof}
Take a point $r\in U_p\cap U_q$ and a relatively compact neighborhood $W$ of $Y_r=\pi^{-1}(r)$ in $V_p\cap V_q$. 
From Lemma \ref{lem:properness_fund}, it follows that there exit a neighborhood $W$ of $r$ in $U_p\cap U_q$ and neighborhood $\Omega$ of $(0, 0, \dots, 0)$ in $\Delta^d$ such that $F_p^{-1}(Z\times W\times \Omega)\subset V_p\cap V_q$. 
By the maximal principle, holomorphic functions $h_q^\nu$'s and $v_q^\mu$'s are constant along a compact submanifolds $F_p^{-1}(Z\times\{(h_p, v_p)\})$ for each $(h_p, v_p)\in W\times \Omega$, from which the assertion holds (Here we have shown the assertion only on $Z\times W\times \Omega$. However it is sufficient by considering the theorem of identity for the holomorphic functions $\partial \Phi^\bullet_{qp}/\partial \zeta$, where $\zeta$ is a coordinate of $Z$).  
\end{proof}

\begin{center}
\includegraphics[width=6cm]{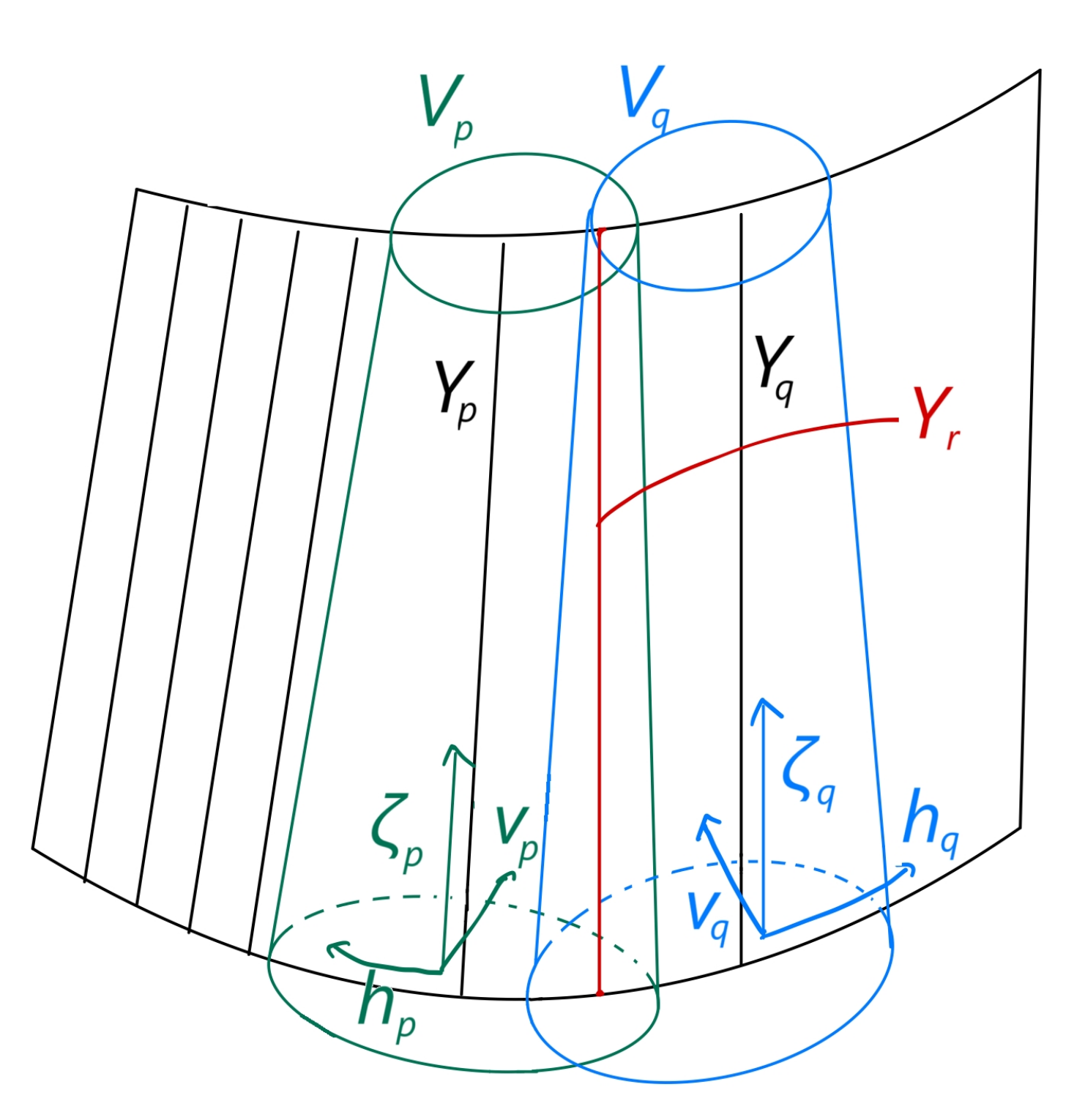}
\end{center}

\subsection{Construction of a new model $(C, M)$}\label{subsection:constr_of_CM}

Assume that $\pi\colon Y\to C$ is as in \S \ref{subsection:someFanoBdls}, and that $Y$ is embedded in a complex manifold $X$ so that the normal bundle $N_{Y/X}$ is a unitary flat vector bundle on $Y$ of rank $d$. 
From Lemma \ref{lem:nice_atlas_coord_change}, it follows that we can construct an another model $C\subset M$ by letting 
\[
M := \left(\bigcup_{p\in C}(U_p\times \Delta^d)\right)/\sim, 
\]
where patching is done by $h_q=\Phi^h_{qp}(h_p, v_p)$ and $v_q=\Phi^v_{qp}(h_p, v_p)$ as in Lemma \ref{lem:nice_atlas_coord_change}. 
By patching the maps
\[
\xymatrix{
V_p \ar[r]^{\cong\ \ \ \ \ \ }_{F_p\ \ \ \ \ \ }& \mathbb{P}^r\times U_p\times \Delta^{d}\ar[r]^{\ \ {\rm Pr}_2} & U_p\times \Delta^{d}, 
}
\]
one has the holomorphic surjective submersion 
\[
\widehat{\pi}\colon V\to M, 
\]
where 
\[
V := \bigcup_{p\in C}V_p. 
\]
Note that 
\[
\widehat{\pi}|_Y = \pi
\]
holds by construction. 
Note that $N_{Y/X}\cong \pi^*N_{C/M}$ holds, since the dual of both hands sides are locally generated by 
\[
\frac{\partial}{\partial v_p^1}, 
\frac{\partial}{\partial v_p^2}, 
\dots, 
\frac{\partial}{\partial v_p^d}. 
\]

As $\widehat{\pi}$ is a $\mathbb{P}^r$-fibration, 
denoting by $M_\bullet$ the open subset $U_\bullet\times \Delta^d\subset M$ and by $\zeta_\bullet$ the fiber coordinates on $V_\bullet$, there exists a holomorphic map $G_{qp}\colon M_p\cap M_q\to {\rm Aut}(\mathbb{P}^r)$ such that the coordinate transitions on $V_p\cap V_q$ are as follows:
\begin{equation}\label{eq:transition_original}
\begin{cases}
h_q=\Phi^h_{qp}(h_p, v_p)\\
v_q=\Phi^v_{qp}(h_p, v_p)\\
\zeta_q = G_{qp}(h_p, v_p)\cdot \zeta_p
\end{cases}.  
\end{equation}

\subsection{Proof of Theorem \ref{cor:main}}
Here we prove Theorem \ref{cor:main} by using the notation in the previous subsections. In what follows we assume that $n=1$, and that $N_{Y/X}$ is Diophantine. 
It follows from \cite[Theorem 1.1]{stolo-gong-tori} that the neighborhood of $C$ in $M$ is full-linearizable since $N_{C/M}$ satisfies the Diophantine condition. Thus, the neighborhoods of $C$ in $M$ can be linearizable. 
Thus we may assume that $\Phi^\bullet_{qp}$'s are as in 
Lemma \ref{lem:nice_atlas_coord_change} are in the form of 
\[
\begin{cases}
h_q=\Phi^h_{qp}(h_p)\\
v_q=\Phi^v_{qp}(v_p) = t_{qp}\cdot v_p
\end{cases}, 
\]
where $t_{qp}$ is an element of ${\rm U}(d)$, or equivalently, we may assume that $M$ is an open neighborhood of the zero section in the total space of the normal bundle $N_{C/M}$. In what follows we denote by $P\colon M \to C$ the holomorphic retraction which is obtained by regarding $M$ as an open subset of $N_{C/M}$ and restricting the projection $N_{C/M}\to C$.

Now the coordinate transitions (\ref{eq:transition_original}) on $V_p\cap V_q$ are improved as follows:
\begin{equation}\label{eq:transition_2}
\begin{cases}
h_q=\Phi^h_{qp}(h_p)\\
v_q=t_{qp}\cdot v_p\\
\zeta_q = G_{qp}(h_p, v_p)\cdot \zeta_p
\end{cases}.  
\end{equation}
For proving the full-linearizablity of a neighborhood of $Y$ in $M$, it is enough to show $G_{qp}(h_p, v_p)=G_{qp}(h_p)$ by modifying $\zeta_p$'s. 
This condition can be reworded as in the statement of Proposition \ref{prop:star_Pr_bdl_over_ellipt_curve} below by using $P^*Y=\{(y, x)\in Y\times M\mid \pi(y)=P(x)\}$. Indeed, as $P(v_p, h_p)=h_p$ holds by using our local coordinates, 
	it follows by definition that the transitions of fiber coordinates of $P^*Y$ never depend on $v_p$'s. 
	More precisely, the coordinate transitions of $P^*Y$ on the preimage of $U_p\cap U_q$ are as follows:
	\[
	\begin{cases}
		h_q=\Phi^h_{qp}(h_p)\\
		v_q=t_{qp}\cdot v_p\\
		\zeta_q = G_{qp}(h_p)\cdot \zeta_p
	\end{cases}, 
	\]
	where $G_{qp}\colon U_p\cap U_q\to {\rm Aut}(\mathbb{P}^r)$ is a holomorphic map, which is nothing but the transition rule of the fiber coordinates of $\pi\colon Y\to C$. 
	Therefore, from Proposition \ref{prop:star_Pr_bdl_over_ellipt_curve}, it follows that we may assume (by changing the coordinates suitably) that the function $G_{qp}(h_p, v_p)$ in (\ref{eq:transition_2}) does not depend on $v_p$.  Therefore the assertion follows.\qed

\begin{proposition}\label{prop:star_Pr_bdl_over_ellipt_curve}
By shrinking $M$ to a small (tubular) neighborhood of the zero section and by replacing $V$ with $\widehat{\pi}^{-1}(M)$, $\widehat{\pi}\colon V\to M$ coincides with $P^*Y\to M$; i.e. there exists a biholomorphism $V\cong P^*Y$ which makes the following diagram commutative:
\[
\xymatrix{
V \ar[rr]^{\cong}\ar[rd]_{\widehat{\pi}}& & P^*Y\ar[ld] \\
& M &
}.
\]
\end{proposition}

Here we need the assumption $n=1$ ($C$ is an elliptic curve) in order to assure that $V$ in the following proof comes from a holomorphic vector bundle. 

\begin{center}
\includegraphics[width=6cm]{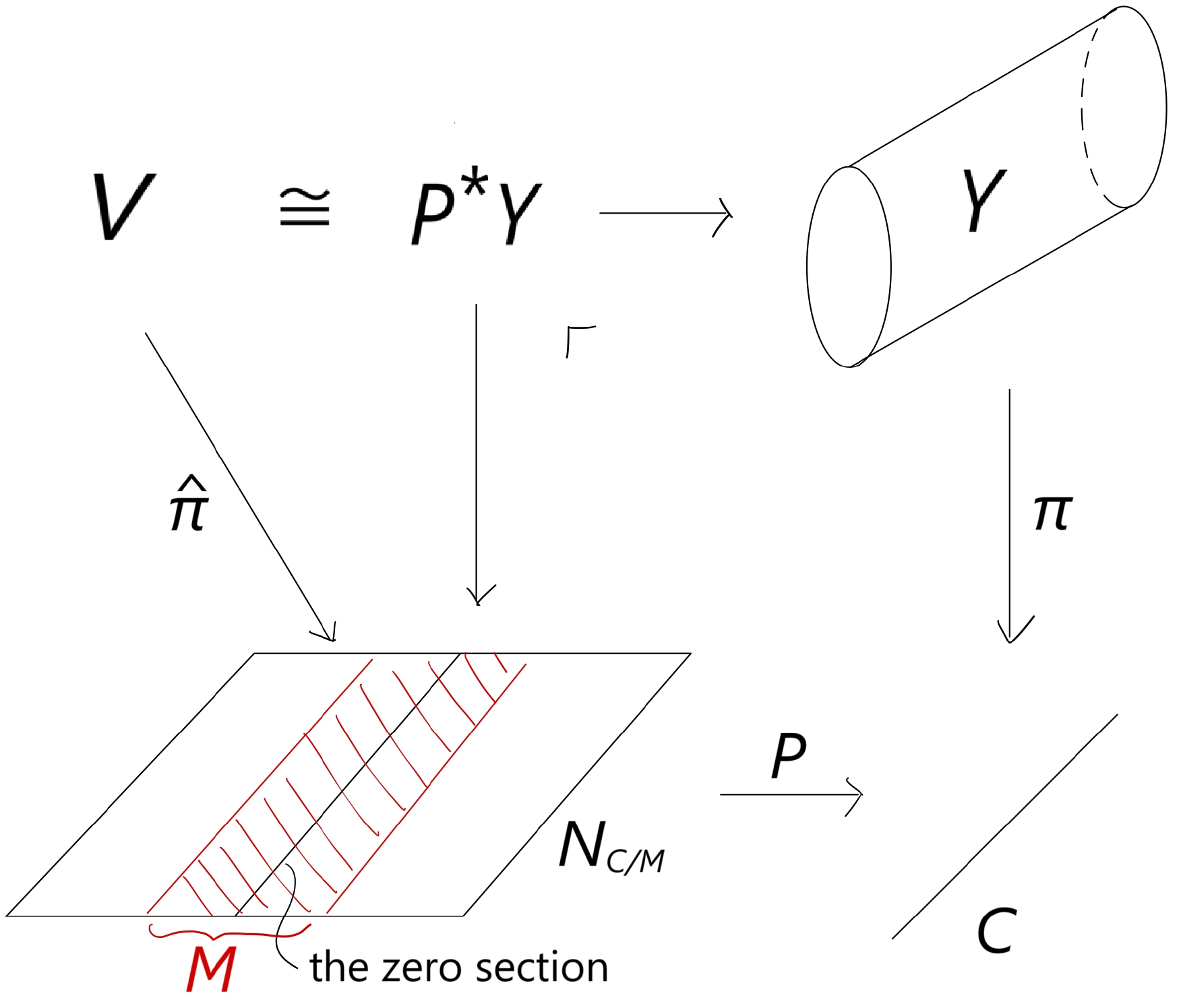}
\end{center}


\begin{proof}[Proof of Proposition \ref{prop:star_Pr_bdl_over_ellipt_curve}]
By shrinking $M$, we will assume that $M$ is a tubular neighborhood of the zero section in what follows. 
By considering the exact sequence 
\[
\cdots \to H^1(M, {\rm GL}(r+1, \mathcal{O}_M)) \to H^1(M, {\rm PGL}(r+1, \mathcal{O}_M)) \to H^2(M, \mathcal{O}_M^*)\to \cdots, 
\]
which is induced from the short exact sequence
\[
0 \to \mathcal{O}_M^* \to {\rm GL}(r+1, \mathcal{O}_M)\to {\rm PGL}(r+1, \mathcal{O}_M)\to 0, 
\]
it follows from the vanishing of $H^2(M, \mathcal{O}_M^*)$ (see Lemma \ref{lem:vanishing_of_H^2O*} below) that there exists a holomorphic vector bundle $F\to M$ such that $V=\mathbb{P}(F)$ (Note that ${\rm GL}(r+1, \mathbb{C})$ is the group for transitions of vector bundles and ${\rm PGL}(r+1, \mathbb{C})={\rm Aut}(Z)$). 

For proving the assertion, it is sufficient to show that the vector bundle $F$ is holomorphically isomorphic to $P^*E$ on a neighborhood of $C$ in $M$, where $E:=F|_C$, which follows from {\bf (Condition $\ast$)}. 
\end{proof}

\begin{lemma}\label{lem:vanishing_of_H^2O*}
Let $M$ be as in the proof of Proposition \ref{prop:star_Pr_bdl_over_ellipt_curve}. 
Then it holds that $H^2(M, \mathcal{O}_M^*)=1$. 
\end{lemma}

\begin{proof}
Consider the exact sequence 
\[
\cdots \to H^2(M, \mathcal{O}_M) \to H^2(M, \mathcal{O}_M^*) \to H^3(M, \mathbb{Z}) \to \cdots, 
\]
which is induced from the exponential exact sequence. As $M$ is homotopic to the elliptic curve $C$, $H^3(M, \mathbb{Z})=0$. 
As $M$ is weakly pseudoconvex and admits a psh exhaustion function whose complex Hessian has ${\rm dim}\,M-1$ positive eigenvalues at each point (consider $\psi:=\|v_j\|^2$, where $v_j=(v_{j, 1}, v_{j, 2}, \dots, v_{j, d})$'s are local defining functions of $C$ in $M$ such that the transitions are unitary), it follows that $H^2(M, \mathcal{O}_M)=H^2(M, K_M\otimes \mathcal{O}_M)=H^{n, 2}(M)=0$ (here we applied \cite[p. 375 Theorem (5.8)]{D})
Thus $H^2(M, \mathcal{O}_M^*)=1$ holds. 
\end{proof}

\begin{proof}[Proof of Theorem \ref{cor_new}]
Let $C, M, P$, and $E$ be as in the proof of Theorem \ref{cor:main}. 
It is sufficient to show that {\bf (Condition $\ast$)} is satisfied when $E$ is holomorphically trivial. 
To show this, take a holomorphic vector bundle $F\to M$ such that $F|_C\cong E$ and 
set $\mathcal{E}:=\mathit{Hom}(F, P^*E)\cong F^*\otimes P^*E$. 
When $E$ is holomorphically trivial, it is clear from the construction that $\mathcal{E}|_C$ is also holomorphically trivial. 
Therefore one can apply Theorem \ref{thm:extension-trivial-tori} to show the existence of a smaller neighborhood $M'$ of $C$ in $M$ such that $\mathcal{E}|_{M'} = \mathit{Hom}(F|_{M'}, (P|_{M'})^*E)$ is holomorphically trivial, from which the assertion follows. 
\end{proof}

\section{Proof of Theorem \ref{thm:extension-trivial-tori}}\label{section:prf_of_thm_extension-trivial-tori}

In this section, we prove Theorem \ref{thm:extension-trivial-tori} by using notation and results from \cite{stolo-gong-tori}.
		 Let $\Lambda$ be a $2n$-dimensional lattice in $\mathbb{C}^n$.
		We may assume that  $\Lambda$ is defined by $2n$  vectors $e_1,\dots, e_n, e_{n+1},\dots,e_{2n}$  of $\mathbb{C}^n$, where $e_i=(0,\ldots, 0,1,0,\ldots, 0)$ with $1$ being at the $i$-th place and $e_{n+i}=(e_{n+i,1},\ldots,e_{n+i,n})$ ,  $1\leq i\leq n$, and the matrix
				$ 
		(\IM{e_{n+i,j}})_{1\leq i,j\leq n}
		$ is invertible 
		The compact complex manifold $C:=\mathbb{C}^n/\Lambda$ is called an ($n$-dimensional) complex torus.
		
		For $\epsilon>0$,  define a (Reinhardt) neighborhood ${\Om}_\epsilon$ of $\overline{{\Om}_0}$ by
		\begin{gather*}
			\om_\epsilon:=\left\{\sum_{j=1}^{2n} t_je_j\colon t\in[0,1)^n\times (-\epsilon,1+\epsilon)^{n}\right\},
			\quad {\Om 
			}_\epsilon:=\{(e^{2\pi i\zeta_1}, \dots e^{2\pi i\zeta_n})\colon\zeta\in \om_\epsilon\}.
		\end{gather*}
		
		With $\Del_r=\{z\in\C\colon|z|<r\}$, we also define
		\beq\label{domains}
		\om_{\epsilon,r}:= \omega_\epsilon\times  \Delta^d_r, 
		\quad \Omega_{\epsilon,r}:=\Om_\epsilon\times  \Delta^d_r. 
		\eeq
		Throughout the paper, a mapping $(z',v')=\psi^0(z,v)$ from $\om_{\e,r}$ into $\C^{n+d}$ that commutes with $z_j\to z_j+1$ for $j=1,\dots, n$ will be identified with a well-defined mapping $(h',v')=\psi(h,v)$ from $\Om_{\e,r}$ into $\C^{n+d}$, where $z,h$ and $z',h'$ are related   as follow in \re{hTj} below.
		\begin{proposition}[{\cite[Proposition 4.3(i)]{stolo-gong-tori}}]\label{transformed-problem}
			Let $C$ be the complex torus and $\pi_{\cyl}\colon\cyl=\C^n/{\Z^n}\to C$ be the covering.
			Let $(M,C)$ be a neighborhood of $C$. Assume that $N_C$ is flat.
			Then one can take $\om_{\e_0,r_0}=\om_{\e_0}\times \Delta_{r_0}^d$ such that $(M,C)$ is biholomorphic to the quotient of $\om_{\e_0,r_0}$ by  $\tau^0_1,\dots,\tau^0_n$. Let $\tau_j$ be the mapping defined on $\Om_{\e_0,r_0}$ corresponding to $\tau_j^0$.
			Then $\tau_1,\dots, \tau_n$ commute pairwise   wherever they are defined, i.e.
			$$
			\tau_i\tau_j(h,v)=\tau_j\tau_i(h,v)\quad  \forall  i\neq j
			$$
			for $(h,v)\in \Om_{\e_0, r_0}\cap \tau_i^{-1}\Om_{\e_0,r_0}\cap \tau_j^{-1}\Om_{\e_0,r_0}$.
		\end{proposition}
The deck transformations of $(\tilde N_C,\cyl)$ are   generated by $n$ biholomorphisms $\hat\tau_1,\dots, \hat\tau_n$ that preserve $\cyl$. Write
\begin{gather}
\hat \tau_j(h,v)=
(T_jh,  M_jv),\quad M_j:={\rm diag}(\mu_{j,1},\dots, \mu_{j,d})\label{tauhat}
\end{gather}
with $h,T_j$ being defined by:
\eq{hTj} h=(e^{2\pi i z_1}, \cdots, e^{2\pi iz_n}), \quad
T_j:= {\rm diag}(\la_{j,1},\dots, \la_{j,n}),\quad \la_{j,k}:=
e^{2\pi ie_{n+j,k}}.
\eeq
		
		A function on $\om_{\e,r}$ that has period $1$ in all $z_j$ is identified with a function on $\Om_{\e,r}$.
		\begin{definition}	
			Set
			$\Om_{\epsilon,r}:=\Om_\epsilon\times\Del_r^d$,
			$
			\tilde \Omega_{\epsilon,r}:=\overline\Omega_{\epsilon,r}\cup\bigcup_{i=1}^n \hat\tau_i(\overline\Omega_{\epsilon,r})$.
			Denote by $\cL A_{\epsilon,r}$ $($resp. $\tilde{\cL A}_{\epsilon, r})$ the set of holomorphic functions on $\overline{\Omega_{\epsilon,r}}$ $($resp. $\overline{\tilde \Omega_{\epsilon,r}})$.
			If $f\in \mathcal{A}_{\epsilon, r}$ $($resp. $\tilde f\in\tilde{\cL A}_{\epsilon, r})$, we set
			\eq{line11}
			\|f\|_{\epsilon,r}:=\sup_{(h,v)\in\Omega_{\epsilon,r}}|  f(h,v)|, \quad |||\tilde f|||_{\epsilon,r}:=\sup_{(h,v)\in\tilde\Omega_{\epsilon,r}}| \tilde f(h,v)|.
			\eeq
		\end{definition}
		Let $C$ be a torus holomorphically embedded in $M$. Under suitable assumption of \cite[Theorem 1.1]{stolo-gong-tori}, a neighborhood of $C$ in $M$ is holomorphically equivalent to a neighborhood of the zero section into its normal bundle $N_{C/M}$.
		
		Recall that
		$$
		\Omega^{ij}_{\epsilon',r'}:= \Omega_{\epsilon',r'}\cap\hat\tau_i^{-1}(\Omega_{\epsilon',r'})\cap\hat\tau_j^{-1}(\Omega_{\epsilon',r'}).
		$$
		Recall that $\|f\|_{\e,r}$ is defined in \re{line11} for a holomorphic function $f\in \cL A_{\e,r}$. For a holomorphic mapping $F=(F_1, F_2, \dots, F_\ell)\in \cL A^\ell_{\e,r}$, define $\|F\|_{\e,r}=\max\{\|F_1\|_{\e,r},\dots, \|F_\ell\|_{\e,r}\}$.
		\begin{definition}\label{def:verticallyDiophantine}
			$N_C$ is said to be unitary flat and {\it vertically Diophantine} if there exist positive constants $D,\tau$ such that  for all $(Q,P)\in \mathbb{N}^d\times \Z^n$,   $|Q|>1$ and all $i=1,\ldots, n$, and  $j=1,\ldots ,d$, we have
			\beq
			\max_{\ell\in \{1,\ldots, n\}}  \left |\la_\ell^P 
			\mu_{\ell}^Q-\mu_{{\ell},j}\right |   >  \frac{D}{(|P|+|Q|)^{\tau}}.\label{dv}
			\eeq
			Furthermore, the $\mu_{i,j}$'s are all of modulus $1$. 
		\end{definition}
		Let us define the operator 
		\begin{align*}
		\hat{\cL L}_i^v:\tilde{\cL A}^{\ell}_{\epsilon, r}&\rightarrow {\cL A}^{\ell}_{\epsilon, r}\\
		G &\mapsto \hat{\cL L}_i^v(G):= G(\hat\tau_i)-G.
		\end{align*} 
		\begin{proposition}[{\cite[Propostion 4.17]{stolo-gong-tori}}]\label{cohomo}
			Assume   $N_C$
			is   
			unitary flat and vertically Diophantine.   Fix $\epsilon_0,r_0,\del_0$  in $(0,1)$. Let $0<\e'<\epsilon<\epsilon_0$, $0<r'<r<r_0$, $0<\del<\del_0$, and $\frac{\del}{\kappa}<\e$. Suppose that
			$F_i=O(|v|^2)\in {\cL A}^{\ell}_{\epsilon, r}$, $i=1,\ldots, n$,
			satisfy
			\beq
			\label{almost-com}  \hat{\cL L}_i^v(F_j)-\hat{\cL L}_j^v(F_i)=0\quad
			\text{on  }  \Omega^{ij}_{\epsilon',r'},\;\text{for all  }i,j.
			\eeq
			There exists a mapping  $G \in  \tilde{\cL A}^{\ell}_{\epsilon
				{-{\del}/{\kappa}},re^{-\delta}}$ 
			such that
			\begin{align}\label{formal2q+1}
				\hat{\cL L}_i^v(G)&=F_i \ \text{on   } \Omega_{\epsilon  {-{\del}/{\kappa}},re^{-\del}}, \;\text{for all  }i.
			\end{align}
			Furthermore, the $G$ satisfies
			\begin{align}\label{estim-sol}
				\|G\|_{\epsilon  {-{\del}/{\kappa}},re^{-\delta}}&\leq  \max_i\|F_{i}\|_{\epsilon,r}\frac{C'}{\delta^{\tau+\nu}},\\
				\label{estim-solcompo}
				\|G\circ\hat\tau_i\|_{\epsilon  {-{\del}/{\kappa}},re^{-\delta}}&\leq  \max_i\|F_{i}\|_{\epsilon,r}\frac{  C'}{\delta^{\tau+\nu}}
			\end{align}
			for some constant 
			$C'$  that is independent of $F,q,\del, r,\e$ and $\nu$ that depends only on $n$ and $\ell$.
			Furthermore, if $F_j(h,v)   
			=J^{2q}F_j(h,v)=O(|v|^{q+1})$ for all $j$, then
			$G$ can be chosen so that
			\eq{}
			G(h,v)=O(|v|^{q+1}), \quad G(h,v)=J^{2q}G(h,v).
			\eeq
		\end{proposition}
		\begin{remark}
		The previous proposition was proved in \cite{stolo-gong-tori} with operator $G\mapsto G(\hat\tau_i)-T_iG$ instead of $\hat{\cL L}_i^v$. Its proof is identical since the numbers to be controlled from below are 
		$$ 
		\left |\la_\ell^P\mu_{\ell}^Q-1\right |=  \left |\la_\ell^P\mu_{\ell}^{Q+e_{j}}-\mu_{{\ell},j}\right |.
		$$
		
		\end{remark}

		\begin{proof}[Proof of Theorem \ref{thm:extension-trivial-tori}]
			Let $C$, $M$ and $\mathcal{E}$ be as in Theorem \ref{thm:extension-trivial-tori}. 
From the discussion above a simple argument implies that, by shrinking $M$ if necessary, there exists a covering map $P_M\colon \widetilde{C}\times \Delta_r^d \to M$ whose deck transformations are generated by $\widehat{\tau}_1, \widehat{\tau}_2, \dots, \widehat{\tau}_n$ such that 
the pull-back $P_M^*\mathcal{E}$ is a holomorphically trivial vector bundle on $\widetilde{C}\times \Delta_r^d$. 

As $P_M^*\mathcal{E}$ is a holomorphically trivial, one can take holomorphic global sections $e_1, e_2, \dots, e_\ell \in \Gamma(\widetilde{C}\times \Delta_r^d, \mathcal{O}_{\widetilde{C}\times \Delta_r^d}(P_M^*\mathcal{E}))$ which holomorphically trivialize $P_M^*\mathcal{E}$, where $\ell$ is the rank of $\mathcal{E}$. 
As the pull-back $\widehat{\tau}_j^*e_\lambda$ is also a holomorphic section of $P_M^*\mathcal{E}$ for each $j\in \{1, 2, \dots, n\}$ and $\lambda\in\{1, 2, \dots, \ell\}$, there exist holomorphic functions $a_{j, \lambda}^\mu\colon \widetilde{C}\times \Delta_r^d\to \mathbb{C}$ such that 
\[
\widehat{\tau}_j^*e_\lambda = \sum_{\mu=1}^\ell a_{j, \lambda}^\mu\cdot e_\mu
\]
holds on $\widetilde{C}\times \Delta_r^d$. 
Let $F_i : \Om_{\epsilon,r}\rightarrow GL(\C^\ell)$ be the 
$\ell\times \ell$ matrix-valued function defined by 
\[ 
F_j := (a_{j, \lambda}^\mu)_{1\leq \lambda, \mu\leq \ell}
\]
for each $j\in \{1, 2, \dots, n\}$. 
Then, as $\widehat{\tau}_i$ and $\widehat{\tau}_j$ are commutative for any for $i, j\in \{1, 2, \dots, n\}$, we have 
\[
\widehat{\tau}_i^*\widehat{\tau}_j^*s = (\widehat{\tau}_j\circ\widehat{\tau}_i)^*s
= (\widehat{\tau}_i\circ\widehat{\tau}_j)^*s = \widehat{\tau}_j^*\widehat{\tau}_i^*s 
\]
for any holomorphic section of $P_M^*\mathcal{E}$. From this, by comparing the coefficients, it follows that $\{F_i\}$ is a collection of $\ell\times \ell$ matrix-valued functions that defines the vector bundle, and that it satisfies the factor of automorphy relations on $\Omega^{ij}_{\epsilon,r}$:
			\beq\label{automorphy}
			F_j(\hat \tau_i)F_i=F_i(\hat \tau_j)F_j,\quad 1\leq i,j\leq n.
			\eeq
			
			In the construction of $e_\lambda$'s above, note that we may assume that $e_\lambda|_{\widetilde{C}\times \{0\}} = \pi_{\widetilde{C}}^*\sigma_\lambda$ for a fixed global holomorphic frame $\sigma_1, \sigma_2, \dots, \sigma_\ell$ of $\mathcal{E}$ on $C$. 
Then, as $\pi_{\widetilde{C}}\circ\widehat{\tau}_j = \pi_{\widetilde{C}}$, we have
\[
\widehat{\tau}_j^*e_\lambda = \widehat{\tau}_j^*\pi_{\widetilde{C}}^*\sigma_\lambda = \pi_{\widetilde{C}}^*\sigma_\lambda = e_\lambda
\]
holds on $\widetilde{C}\times\{0\}$. 
Thus we have $F_i(h,0)=I$ and $F_i(h,v)=I+f_i(h,v)$. Let us prove that we can find $F\colon\Om_{\epsilon',r'}\rightarrow GL(\C^\ell)$, with $0<\epsilon'\leq \epsilon$ and $0<r'<r$ such that $\Phi(h,0)=I$, $\Phi=I+\phi$, $\phi$ is matrix valued function the coefficients of which belongs to $\tilde{\cL A}_{\epsilon',r'}$ and on $\Omega^{i}_{\epsilon',r'}:=\Omega_{\epsilon',r'}\cap\hat\tau_i^{-1}(\Omega_{\epsilon',r'})$, we have
			\beq\label{trivial}
			\Phi(\hat\tau_i(h,v))F_i(h,v)\Phi(h,v)^{-1}=I,\quad i=1,\ldots, n.
			\eeq
			
			We will solve \re{trivial} by iterations thought a Newton scheme. Let us set
			$$
			\del_{k}:=
			\frac{\del_0}{(k+1)^2}, 
			\quad r_{k+1}:=r_ke^{-5\del_k}, \quad \epsilon_{k+1}:= \epsilon_k
			{ {-\frac{5\del_k}{\kappa}}},\;m=2^k-1.
			$$
			We have $\epsilon_{k}>\frac{\epsilon_0}{2}$ and $r_k>\frac{r_0}{2}$, $k\geq 0$ (see \cite[\S 4.3]{stolo-gong-tori}).
			
			We assume that $F_{i,m+1}=I+f_{i,m+1}$ with $f_i=O(v^{m+1})$ and $\|f_{m+1}\|_{\epsilon,r}:=\max_i\|f_{i,m+1}\|_{\epsilon,r}\leq \delta_k^{\mu}$ for some $\mu$ to be fixed below.
			Let us find $\Phi_{m+1}:=I+\phi_{m+1}$ such that $\phi_{m+1}=O(v^{m+1})\in\tilde{\cL A}_{\epsilon',r'}$, for some $0<\epsilon'<\epsilon$, $0<r'<r$, satisfying, on $\Omega^{i}_{\epsilon',r'}$,
			$$
				\Phi_{m+1}(\tau_i(h,v))F_{i,m+1}(h,v)\Phi_{m+1}(h,v)^{-1}=I+f_{i,m+1}^+,\quad i=1,\ldots, n.
			$$
			where $f_{i,m+1}^+=O(v^{2m+1})\in {\cL A}_{\epsilon',r'}^{\ell^2}$. Expanding the previous equality, we have
				\begin{align}
				\hat{\cL L}_i^v(\phi_{m+1})=&f_{i,m+1}^+-f_{i,m+1}\label{conj-equ}\\
				&+f_{i,m+1}^+\phi_{m+1}-\phi_{m+1}(\hat\tau_i)f_{i,m+1}=:R_{i,{m+1}}\nonumber
			\end{align}
			Let us set $\tilde R_i:=J^{2m}(R_{i,{m+1}})$. We have $\tilde R_i=-J^{2m}(f_{i,m+1})$. We recall that right composition by any $\hat \tau_i$ preserve the degree w.r.t. $v$. Hence, $J^{2m}(g(\hat \tau_i))=J^{2m}(g)(\hat\tau_i)$ for any appropriate function $g$.  Using \re{automorphy}, we obtain
			\begin{eqnarray*}
			\hat{\cL L}_i^v(J^{2m}(R_{j,{m+1}}))&=&J^{2m}(-f_{j,m+1}(\hat \tau_i)+f_{j,m+1})\\
			&=&J^{2m}(-F_{j,m+1}(\hat \tau_i)+F_{j,m+1})\\
			&=&J^{2m}(-F_{j,m+1}(\hat \tau_i)F_iF_i^{-1}+F_{j,m+1})\\
			&=&J^{2m}(-F_{i,m+1}(\hat \tau_j)F_jF_i^{-1}+F_{j,m+1})\\
\end{eqnarray*}
As $F_jF_i^{-1}= I+f_j-f_i+O(v^{2m+2})$, we have $F_{i,m+1}(\hat \tau_j)F_jF_i^{-1}= I +f_{i,m+1}(\hat\tau_j)+ f_j-f_i+O(v^{2m+2})$. Therefore, the previous computations show that
\begin{align*}
\hat{\cL L}_i^v(J^{2m}(R_{j,{m+1}}))&=J^{2m}(-F_{i,m+1}(\hat \tau_j)F_jF_i^{-1}+F_{j,m+1})\\
&= J^{2m}(-f_{i,m+1}(\hat \tau_j)+f_{i,m+1})= \hat{\cL L}_j^v(J^{2m}(R_{i,{m+1}})).
\end{align*}
According to \rp{cohomo}, there exists $\phi_{m+1}\in \tilde{\cL A}^{\ell}_{\epsilon{-{\del}/{\kappa}},re^{-\delta}}$ such that
\begin{align}
\hat{\cL L}_i^v(\phi_{m+1})=\tilde R_i:=&-J^{2m}(f_{i,{m+1}}),\;\text{ for all  }1\leq i\leq n,\nonumber\\ 
	\|\phi_{m+1}\|_{\epsilon  {-{\del}/{\kappa}},re^{-\delta}}, &\|\phi_{m+1}\circ\hat\tau_i\|_{\epsilon  {-{\del}/{\kappa}},re^{-\delta}}\leq  \max_i\|f_{i,m+1}\|_{\epsilon,r}\frac{C'}{\delta^{\tau+\nu}},\label{estim-cohom}\\
	\phi_{m+1}=O(v^{m+1}),&\quad  \phi_{m+1}=J^{2m}(\phi_{m+1}).\nonumber
\end{align}
According to \re{conj-equ}, we then have on $\Om_{\epsilon{-{\del}/{\kappa}},re^{-\delta}}$,
\beq\label{remainder}
f_{i,m+1}^+:=\left(\phi_{m+1}(\hat\tau_i)f_{i,m+1}+(f_{i,m+1}-J^{2m}(f_{i,m+1}))\right)(I+\phi_{m+1})^{-1}.
\eeq
Recalling that $\|f_{m+1}\|_{\epsilon,r}\leq \delta_k^{\mu}$ with $m=2^k-1$.
Hence, using Schwarz inequality, we have for $0<\delta< \kappa\epsilon/2$,
\begin{align*}
\|f_{m+1}^+\|_{\epsilon {-{2\del}/{\kappa}},re^{-2\delta}}&\leq\left(\frac{C'\|f_{m+1}\|_{\epsilon {-{\del}/{\kappa}},re^{-\delta}}^2}{\delta^{\tau+\nu}}+e^{-(2m+1)\delta}\|f_{m+1}\|_{\epsilon {-{\del}/{\kappa}},re^{-\delta}}\right)\frac{1}{1-\left(\frac{C'\|f_{m+1}\|_{\epsilon {-{\del}/{\kappa}},re^{-\delta}}}{\delta^{\tau+\nu}}\right)}\\
&\leq \left(\frac{C'e^{-2(m+1)\delta}\|f_{m+1}\|^2_{\epsilon,r}}{\delta^{\tau+\nu}}+e^{-(3m+2)\delta}\|f_{m+1}\|_{\epsilon ,r}\right)\frac{1}{1-\left(\frac{C'e^{-(m+1)\delta}\|f_{m+1}\|_{\epsilon,r}}{\delta^{\tau+\nu}}\right)}\\
&\leq \left(C'e^{-2(m+1)\delta}\delta_k^{2\mu}\delta^{-\tau-\nu}+e^{-(3m+2)\delta}\delta_k^{\mu}\right)\frac{1}{1-\left(C'e^{-(m+1)\delta}\delta_k^{\mu}\delta^{-\tau-\nu}\right)}.\\
\end{align*}
Let us set $\delta:=\delta_k$. We can assume that $\mu$ is large enough so that, for all $k\geq 0$, $C'e^{-(m+1)\delta_k}\delta_k^{\mu-\tau-\nu}<\frac{C'}{2^k}\delta_k^{\mu-\tau-\nu-1}<\frac{1}{2}$.
Let us check that 
$$
C'e^{-2(m+1)\delta_k}\delta_k^{2\mu}\delta_k^{-\tau-\nu}+e^{-(3m+2)\delta_k}\delta_k^{\mu}< e^{-(m+1)\delta_k}\delta_k^{\mu}\leq \frac{1}{2}\delta_{k+1}^{\mu}.
$$
Indeed, we have
$$
e^{-(m+1)\delta_k}\leq \frac{(k+1)^2}{2^k\delta_0}\leq \frac{1}{2^{2\mu+1}}\leq \frac{1}{2}\left(\frac{1}{1+\frac{1}{k+1}}\right)^{2\mu}=\frac{1}{2}\left(\frac{\delta_{k+1}}{\delta_k}\right)^{\mu},
$$
if $k\geq k_0$ is large enough.
Therefore we have, on $\Om_{\epsilon_{k+1},r_{k+1}}^i$
$$
\Phi_{m+1}(\tau_i(h,v))F_{i,m+1}(h,v)\Phi_{m+1}(h,v)^{-1}=F_{i,2m+1}=:I+f_{i,2m+1},\quad i=1,\ldots, n.
$$
where 
$$
f_{i,2m+1}:=f^+_{i,m+1}=O(v^{2m+1}),\quad \|f_{2m+1}\|_{\epsilon_{k+1},r_{k+1}}\leq \delta_{k+1}^{\mu},\;\;m=2^k-1.
$$
Furthermore, according to \re{estim-cohom}, we have
\begin{align*}
	\|\phi_{m+1}\|_{\epsilon_k  {-{2\del_k}/{\kappa}},r_ke^{-2\delta_k}}
	&\leq  e^{-2(m+1)\delta_k}\|f_{m+1}\|_{\epsilon_k  ,r_k}\frac{C'}{(2\delta_k)^{\tau+\nu}},\\
	&\leq \frac{C'}{2^{\tau+\nu}}e^{-2(m+1)\delta_k}\delta_k^{\mu-\tau-\nu}\leq \frac{1}{2^{\tau+\nu+1}}e^{-(m+1)\delta_k}
\end{align*}
The same inequality holds for $\|\phi_{m+1}\circ\hat\tau_i\|_{\epsilon_k  {-{2\del_k}/{\kappa}},r_ke^{-2\delta_k}}$.
As the sum $\sum_{k\geq 0}e^{-\frac{\delta_02^k}{(k+1)^2}}$ converges, the product $\Phi_{\infty}:=\cdots (I+\phi_k)(I+\phi_{k-1})\cdots(I+\phi_1)=:I+\phi_{\infty}$ converges uniformly on $\tilde{\cL A}_{\e_0/2,r_0/2}$ and on $\Om_{\e_0/2,r_0/2}^i$,  we have 
$$
\Phi_{\infty}(\tau_i(h,v))F_{i}(h,v)\Phi_{\infty}(h,v)^{-1}=I,\quad i=1,\ldots, n,
$$
which means that the vector bundle ${\mathcal E}\rightarrow M$ is holomorphically trivial by shrinking $M$.
\end{proof}


\bigskip

\begin{thebibliography}{99}
	
\bibitem{arnold-elliptic}
V.~I. Arnol'd.
\newblock Bifurcations of invariant manifolds of differential equations, and
normal forms of neighborhoods of elliptic curves.
\newblock {\it Funkcional. Anal. i Prilo\v{z}en.}, 10(4):1--12, 1976.

\bibitem[A]{A} V. I. Arnold, Geometrical Methods in the Theory of Ordinary Differential Equations, Grundlehren der mathematischen Wissenschaften (GL, volume 250). 



\bibitem[D]{D} J.-P. Demailly, Complex analytic and differential geometry, \url{http://www-fourier.ujf-grenoble.fr/~demailly/manuscripts/agbook.pdf}.

\bibitem[GS21]{GS21} X. Gong, L. Stolovitch, Equivalence of Neighborhoods of Embedded Compact Complex Manifolds and Higher Codimension Foliations, Arnold Math J. {\bf 8} (2022), 61--145. 

\bibitem[GS]{stolo-gong-tori} X. Gong, L. Stolovitch, On neighborhoods of embedded complex tori, Math. Ann. (2024), \url{https://doi.org/10.1007/s00208-024-02975-w}. 

\bibitem[Gr]{grauert-neg}
H.~Grauert.
\newblock \"{U}ber {M}odifikationen und exzeptionelle analytische {M}engen.
\newblock {\it Math. Ann.}, 146:331--368, 1962.

\bibitem[G]{G} P. Griffiths, The Extension Problem in Complex Analysis II; Embeddings with Positive Normal Bundle, Amer. J. Math. {\bf 88}, 2 (1966), 366--446.

\bibitem[HR]{hironaka-rossi}
H.~Hironaka and H.~Rossi.
\newblock On the equivalence of imbeddings of exceptional complex spaces.
\newblock {\it Math. Ann.}, 156:313--333, 1964.

\bibitem[H]{H} D. Huybrechts, Lectures on K3 surfaces, Cambridge Studies in Advanced Mathematics, {\bf 158}., Cambridge University Press, Cambridge (2016).

\bibitem[I]{iena-05}
Oleksandr Iena.
\newblock Vector bundles on elliptic curves and factors of automorphy.
\newblock {\em Rend. Istit. Mat. Univ. Trieste}, 43:61--94, 2011.

\bibitem[IP]{ilyashenko-tori}
Y.~S. Ilyashenko and A.~S. Pjartli.
\newblock Neighborhoods of zero type imbeddings of complex tori.
\newblock {\it Trudy Sem. Petrovsk.}, (5):85--95, 1979.

\bibitem[KO]{KO} S. Kobayashi, T. Ochiai, On complex manifolds with positive tangent bundles, J. Math. Soc. Japan, {\bf 22} (4) (1970), 499--525. 

\bibitem[KNS]{KNS} K. Kodaira, L. Nirenberg, D.~C. Spencer, On the Existence of Deformations of Complex Analytic Structures, Ann. of Math., {\bf 68}, 2 (1958), 450--459. 


\bibitem[K]{Koi} T. Koike, Higher codimensional Ueda theory for a compact submanifold with unitary flat normal bundle, Nagoya Math. J., {\bf 238} (2020), 104--136. 

\bibitem[K2]{Koi2} T. Koike, $\overline{\partial}$ cohomology of the complement of a semi-positive anticanonical divisor of a compact surface, Math. Z., {\bf 308}, 22 (2024). 

\bibitem[KU]{KU} T. Koike, T. Uehara, A gluing construction of projective K3 surfaces, Epijournal de Geometrie Algebrique, 6 juillet 2022, Volume {\bf 6} - \url{https://doi.org/10.46298/epiga.2022.volume6.8504}. 


\bibitem[MK]{KM} J. Morrow, K. Kodaira, Complex Manifolds, AMS Chelsea Publishing (1971). 

 

\bibitem[S]{S} Y. T. Siu, Every Stein subvariety admits a Stein neighborhood. Invent. Math. {\bf 38} (1976/77), 89--100.

\bibitem[SW]{stolo-wu} L. Stolovitch and X. Wu, Ueda foliation problem for complex tori. arXiv 2403.17682 (2024)-\url{https://arxiv.org/abs/2403.17682}.


\bibitem[T]{T} H. Tsuji, Logarithmic Fano Manifolds Are Simply Connected, Tokyo J. Math. {\bf 11}(2): 359--362 (1988). 

\bibitem[U]{U} T. Ueda, On the neighborhood of a compact complex curve
with topologically trivial normal bundle, J. Math. Kyoto Univ. {\bf 22} (1982/83), 4, 583-607. 

\end{thebibliography}
\end{document}